\documentclass[]{amsart}
\usepackage[]{amsmath, amsthm, amsfonts }

\newcommand {\C} {{\mathbb C}}

\newcommand {\PP} {{\mathbb P}}
\newcommand {\CC} {{\mathcal C}}
\newcommand {\br} {\overline}
\newcommand {\s} {{\mathcal S}}

\newtheorem{thm}{Theorem}
\newtheorem{cor}{Corollary}
\newtheorem{lemma}{Lemma}
\newtheorem{prop}{Proposition}

\begin{document}

\title{
Kodaira dimension of symmetric powers
}
\author{
        Donu Arapura    
}
\author{
Sviatoslav Archava
}
\address{Department of Mathematics\\
Purdue University\\
West Lafayette, IN 47907\\
U.S.A.}
\thanks{First author partially supported by the NSF}

\maketitle

We work over the complex numbers.
When $X$ is a smooth projective curve of genus
$g$, elementary arguments show that the $d$th
symmetric power $S^dX$ is uniruled as soon
as $d>g$, and therefore that  the plurigenera  vanish.
When the dimension of $X$ is greater
than one, the situation is quite different.
Using ideas of Roitman \cite{roitman,
roitman2} and Reid \cite{reid}, we prove:

\begin{thm} Let $X$ be smooth projective
variety with $n = dim X> 1$.
Let $\Sigma_d$ be a desingularization
of $S^dX$, then there are isomorphisms
$$S^dH^0(X, \omega_X^{\otimes m})\cong 
H^0(\Sigma_d,\omega_{\Sigma_d}^{\otimes m})$$
whenever $mn$ is even.
\end{thm}

\begin{cor}
With the previous assumptions, the
$m$th plurigenus
$$ P_m(\Sigma_d) = 
\left(
\begin{array}{c}d + P_m(X)-1 \\ d\\ \end{array}
\right)
$$
whenever $mn$ is even.
 The Kodaira dimension
$\kappa(\Sigma_d) = d\kappa(X)$.
\end{cor}

\begin{proof}
The first formula is an immediate consequence of the
theorem. It implies that 
$$P_m(\Sigma_d) = O(P_m(X)^d) = O(m^{d\kappa(X)})$$
which yields the second formula.
\end{proof}

Recall that a projective variety $Z$ is uniruled provided
there exists a variety $Z'$ and dominant rational map 
$Z'\times \PP^1\dashrightarrow Z$. The reference
\cite{kollar} is  more than adequate for standard
properties of uniruled varieties.

\begin{cor}
If $X$ has nonnegative Kodaira dimension
then $S^dX$ is not uniruled for any $d$.
\end{cor}

\begin{proof}
 Since uniruledness is a birational property,
 it is enough to observe that $\Sigma_{d}$ is not
 uniruled because it has nonnegative Kodaira
dimension.
\end{proof}

The most interesting corollaries involve genus
estimates for curves lying on $X$. The phrase ``$d$ general points
of $X$ lie on an irreducible curve with
genus $g$ normalization''
 will mean that there is an irreducible quasiprojective
family $\CC\to T$ of 
smooth projective genus $g$ curves and a morphism $\CC\to X$ which is
a generically one to one on the fibers $\CC_t$ and such that the
morphism from the 
relative symmetric power 
$${\mathcal S}^d\CC := \CC\times_T\CC\times_T\ldots \CC/S_d $$
to $S^dX$ is dominant.

\begin{cor}
Suppose that the Kodaira dimension of $X$ is
nonnegative and that  $d$ general points
lie on  an irreducible curve with genus $g$ normalization, then
$g \ge d$.
\end{cor}

\begin{proof}
Assume the contrary that $g < d$, and let $\CC\to T$ be
the corresponding family. Then each fiber $\s^d\CC_t$
is a projective space bundle over the Jacobian $J(\CC_t)$ 
by Abel-Jacobi; in particular, it's uniruled. Therefore
$\s^d\CC$ and hence $S^dX$ are uniruled, but this
contradicts the previous corollary.
\end{proof}

\begin{cor}
Suppose that  $X$ has general type  and that  $d$ general points
lie on  an irreducible curve with genus $g$ normalization, then
$g > d$.
\end{cor}

\begin{proof}
We assume that $g\le d$ for some family $\CC\to T$.
By the previous corollary, we may
suppose that $g=d$. Denote the maps $\s^d\CC\to S^dX$
and  $\mathcal{S}^d\CC\to T$ by $p$ and $\pi$ respectively.
The map $p$ is dominant and generically injective on the fibers of 
$\pi$.
If a general fiber $p^{-1}(Z)$ has positive dimension, then
an irreducible hyperplane section $H\subset T$
 meets $\pi(p^{-1}(Z))$.
 Therefore $\mathcal{S}^d\CC\times_T H\to S^dX$
is still dominant, and we may replace $T$ by $H$ and
$\mathcal{S}^d\CC$ by the fiber product. By continuing in this way, we
can assume that $p$ is generically finite. Choose
a desingularization $\br T$ of a compactification of $T$,
and a nonsingular compactification of $S$ of
$\mathcal{S}^d\CC\times_T {\br T}$ such that $p$ extends to a morphism
of $S$ to a desingularization $\Sigma_d$ of $S^dX$.
We then have $\kappa(\Sigma_d)\le \kappa(S)$ which implies that
$S$ has general type.
On the other hand the general fiber of $S\to \br T$ is
$\mathcal{S}^d\CC_t$ is birational  to an Abelian variety by
the theorems of Abel and Jacobi. This implies that 
$$\kappa(S)\le dim\br T+\kappa(\mathcal{S}^d\CC_t)= dim \br T < dim S$$ 
by \cite[2.3]{mori},
but this  is impossible since $S$ has general type.
\end{proof}

\section{Proof of main theorem}

Recall \cite{reid} that a 
  variety $Y$ has  canonical singularities provided that 

\begin{enumerate}

\item $Y$ is normal

\item  $\omega_Y^{[r]} := (\omega_Y^{\otimes r})^{**}$ is
locally free for some $r>0$, where 
$\omega_Y = (\Omega_Y^{dimY})^{**}$.

\item If $f:Y'\to Y$ is a resolution
of singularities, then $f_*\omega_{Y'}^{\otimes r} = \omega_Y^{[r]}$.

\end{enumerate}
The smallest such $r$ is called the index. 
If $Y$ is canonical, then the  third condition holds for all
$r\ge 1$ \cite[1.3]{reid}, thus the index is the smallest $r$ for which
$\omega_Y^{[r]}$ is locally free.
It's enough to test the last  condition for a
particular resolution of singularities. This
 condition is equivalent to a more widely used
condition involving
pullbacks of canonical divisors.

\begin{lemma}
Let $Z$ be a smooth variety on which a finite group $G$
acts. Let $Y = Z/G$. If $Y$ has canonical singularities
of index dividing $r$, then 
$$H^0(Y,\omega_Y^{[r]})\subseteq H^0(Z, \omega_X^{\otimes r})^G.$$
If the fix point locus has codimension greater than
one, equality holds.
\end{lemma}

\begin{proof}
Construct a commutative diagram
$$\begin{array}{ccc}
Z'        & \to & Z\\
g'\downarrow&     & g\downarrow \\
Y'        & \stackrel{f}{\longrightarrow} & Y\\
\end{array} $$
where $Y'\to Y$ is a desingularization, and $Z'$ is
a $G$-equivariant desingularization of the fiber product.
Then there are inclusions
$$\omega_Y^{[r]}= f_*\omega_{Y'}^{\otimes r} \subseteq  
(f\circ g')_*(\omega_{Z'}^{\otimes r})^G =
g_*(\omega_{Z}^{\otimes r})^G. $$
This implies the first part of the lemma. The
inclusion $\omega_Y^{[r]} \subseteq  
g_*(\omega_{Z}^{\otimes r})^G$ is an equality  on the complement of 
the fixed point locus. Since $\omega_Y^{[r]}$ is locally
free (hence reflexive) and $g_*(\omega_{Z}^{\otimes r})^G$
is torsion free, the second statement follows.
\end{proof}

\begin{prop}
Let $X$ be a smooth variety of dimension $n > 1$,
then $S^dX$ has canonical singularities of index one
if $n$ is even, and canonical singularities of
index at most $2$ if $n$ is odd.
\end{prop}

\begin{proof}
As the result is local analytic for $X$, we
may  replace it by $\C^n$. Consider the
action of the symmetric group $S_d$ on
$\C^{nd}= \C^n\times \ldots \times\C^n$ by permutation
of factors, and let $h:S_d\to  GL_{nd}(\C)$
be the corresponding homomorphism. This action is equivalent to a
direct sum of $n$ copies of the standard
representation $\C^d$ where $S_d$ acts via permutation matrices. Therefore 
$h(S_d)$ does not contain any quasi-reflections (because $n>1$)
and $det(h(\sigma)) = sign(\sigma)^n$.
When $n$ is even, $h(S_d)\subset SL_{nd}(\C)$.
This implies that $S^dX= \C^{nd}/S_d$ is
Gorenstein by  \cite{watanabe}, and therefore
canonical of index one \cite[1.8]{reid}.

The case when $n$ is odd is more laborious. 
For any element $\sigma \in S_{d}$ of order $r$,
define $S(\sigma)$ as follows: choose a primitive
$r$th root of unity $\epsilon$ and
express the eigenvalues 
of  $h(\sigma)$ as $\lambda_i = \epsilon^{a_i}$
where $0\le a_i< r$, set $S(\sigma) = \sum\, a_i$.
By
\cite[3.1]{reid}, to prove that $S^dX$ is canonical
it will suffice to verify that
$S(\sigma) \ge r$ for every element $\sigma$ of order $r$.
Let $\C^d$ be the permutation representation  of $S_d$.
If $e_1,\ldots e_d$ is the standard basis then
$\sigma\cdot e_i = e_{\sigma(i)}$. If 
$\epsilon$ is a primitive $r$th root of
unity, then it is easy to see that the 
eigenvectors of the cycle $\sigma = (12\ldots r)$ acting on $\C^d$
are $e_1+\epsilon^{i}e_2+\ldots\epsilon^{(r-1)i}e_{r}$
and $e_{r+1},\ldots e_{d}$. Therefore the
nonunit eigenvalues
are $\epsilon,\ldots, \epsilon^{r-1}$ and these occur
with multiplicity one. Hence $S(\sigma) = nr(r-1)/2 \ge r$
as required. The general case is similar. Let $\sigma$ be a permutation
of order $r$ and $\epsilon$ as before. Write $\sigma$
as a product of disjoint cycles of length $r_i$.
Therefore $r$ is the least common multiple of the $r_i$,
and let $r_i' = r/r_i$. A list (with possible repetitions)
of the  nonunit eigenvalues of
$\sigma$ acting on $\C^d$ 
is 
$$\epsilon^{r_1'},\ldots\epsilon^{r_1'(r_1-1)},\epsilon^{r_2'}\ldots
\epsilon^{r_2'(r_2-1)},\ldots$$
Therefore 
$$S(\sigma) = \frac{n}{2}[r_1'r_1(r_1-1) + r_2'r_2(r_2-1)+\ldots] \ge r$$
and this proves that $S^dX$ is canonical. It remains to check that
the index is at most $2$. For this it suffices to observe that if $x_i$ are
coordinates on $\C^{nd}$, then
$$(dx_1\wedge\ldots \wedge dx_{nd})^{\otimes 2}$$
is $S_d$ invariant. This determines a generator
of $\omega_{S^dX}^{[2]}(X)$, which shows that this module
is free.

\end{proof}
 
\begin{proof}[Proof of main theorem]
Let  $m$ be  an integer such that $mn$ is even (hence
a multiple of the index of $S^dX$).
Then 
$$H^0(\omega_{\Sigma_d}^{\otimes m})=H^0(\omega_{S^dX}^{[m]}) 
=H^0(\omega_{X^d}^{\otimes m})^{S_d} $$
By K\"unneth's formula, this equals
$$[H^0(\omega_X^{\otimes m})\otimes \ldots \otimes H^0(\omega_X^{\otimes 
m})]^{S_d}
= S^dH^0(\omega_X^{\otimes m}).$$
\end{proof}



\begin{thebibliography}{ABC}
\bibitem[K]{kollar} J. Koll\'ar,
{\em Rational curves on algebraic varieties}
Springer-Verlag (1996)

\bibitem[Mo]{mori} S. Mori,
{\em Classification of higher dimensional varieties}
Algebraic Geometry Bodwoin 1985, AMS (1987)

\bibitem[M]{mumford} D. Mumford, 
{\em Rational equivalence of $0$-cycles on surfaces}
J. Math. Kyoto Univ. 9 (1969)

\bibitem[Re]{reid} M. Reid,
{\em Canonical $3$-folds}, Journe\'es
G\'eometrie Alg\'ebrique, Sitjoff and Noordorff (1980)

\bibitem[Ro1]{roitman} A. A. Roitman
{\em On $\Gamma$-equivalence of zero-dimensional cycles}
USSR Sbornik 15 (1971)

\bibitem[Ro2]{roitman2} A. A. Roitman
{\em Rational equivalence of $0$-dimensional cycles}
USSR Sbornik 18 (1972)


\bibitem[W]{watanabe} K. Watanabe,
{\em Certain invariant subrings are Gorenstein II}
Osaka Math. J. 11 (1974)
\end{thebibliography}
\end{document}